\documentclass[11pt]{article}
\usepackage{amsmath,amssymb,graphicx,float}

\numberwithin{equation}{section}
\def\eqn{Equation~}
\def\eqns{Equations~}

\begin{document}

\title{\bf The Imani Periodic Functions:\\
Genesis and Preliminary Results}

\author{Ronald E. Mickens\vspace{6pt}\thanks{Email: 
rmickens@cau.edu}\\\vspace{6pt}
Department of Physics, Clark Atlanta University, Atlanta, GA 30314, USA}
\maketitle

\begin{abstract}
The Leah-Hamiltonian, $H(x,y)=y^2/2+3x^{4/3}/4$, is introduced as a
functional equation for $x(t)$ and $y(t)$. By means of a nonlinear 
transformation to new independent variables, we show that this
functional equation has a special class of periodic
solutions which we designate the Imani functions.  The 
explicit construction of these functions is done such that they 
possess many of the general properties of the standard trigonometric
cosine and sine functions.  We conclude by providing a listing of a number 
of currently unresolved issues relating to the Imani functions.
\bigskip



\end{abstract}

\section{Introduction}\label{sec1}

The Leah differential equation is \cite{9,10}
\begin{equation}\label{eq1.1}
\frac{d^2x}{dt^2}+x^{1/3}=0.
\end{equation}
This equation models a ``truly nonlinear oscillator" and arises in the
study of the dynamics of nonlinear oscillatory systems \cite{5}.
The first-integral or Hamiltonian corresponding to \eqn\eqref{eq1.1} 
is \cite{4,5}
\begin{equation}\label{eq1.2}
H(x,y)=\frac{y^2}2+\left(\frac34\right)x^{4/3}=\frac34,
\end{equation}
where
\begin{equation}\label{eq1.3}
y(t)=\frac{dx(t)}{dt},
\end{equation}
and the following initial conditions are selected
\begin{equation}\label{eq1.4}
x(0)=1,\quad y(0)=0.
\end{equation}
The equations of the trajectories, $y=y(x)$, in the $(x,y)$ phase-space
are determined by the first-order differential equation \cite{4,5}
\begin{equation}\label{eq1.5}
\frac{dy}{dx} =-\frac{x^{1/3}}y.
\end{equation}

Using standard techniques from the qualitative theory of differential
equations, it can be easily shown that the following statements are true:

(i) Both $x(t)$ and $y(t)$ are bounded and periodic, with a period $T$
which can be explicitly calculated \cite{5,9}.

(ii) $x(t)$ and $y(t)$ are, respectively, even and odd, i.e.,
\begin{equation}\label{eq1.6}
x(-t)=x(t),\quad y(-t)=-y(t).
\end{equation}

(iii) The Taylor series for $x(t)$ only exists over the interval
\begin{equation}\label{eq1.7}
-\frac T4 < t < \frac T4.
\end{equation}

(iv) The Fourier series for $x(t)$ and $y(t)$ contain only odd 
harmonics and have the representations
\begin{subequations}\label{eq1.8}
\begin{align}
x(t)&=\sum^\infty_{k=0} a_k\cos(2k+1)\left(\frac{2\pi}T\right)t,
\label{eq1.8a}\\
y(t)&=-\sum^\infty_{k=0} (2k+1)\left(\frac{2\pi}T\right)a_k
\sin(2k+1)\left(\frac{2\pi}T\right)t.
\label{eq1.8b}\end{align}
\end{subequations}

\centerline{\hbox to 3in{\hrulefill}}

Note that since the Hamiltonian, $H(x,y)$, see \eqn\eqref{eq1.2}, is
constant, then upon taking its time derivative, it follows that
\begin{align}\label{eq1.9}
\frac{dH(x,y)}{dt} &= \frac{\partial H}{\partial x}\frac{dx}{dt}
+\frac{\partial H}{\partial y}\frac{dy}{dt}\nonumber\\
&=x^{1/3}\frac{dx}{dt} + y\frac{dy}{dt}\nonumber\\
&=x^{1/3}y + y\frac{d^2x}{dt^2}\nonumber\\
&=y\left(x^{1/3}+\frac{d^2x}{dt^2}\right)=0.
\end{align}
Consequently, this provides a derivation of the equation of motion, as
indicated in \eqn\eqref{eq1.1}.

However, in general, the Hamiltonian does not provide a unique set of
equations of motion.  This can be seen easily by considering the first
line of the expression given in \eqn\eqref{eq1.9}, i.e., for
\begin{equation}\label{eq1.10}
H(x,y)=\mbox{constant},
\end{equation}
then
\begin{equation}\label{eq1.11}
\frac{dH}{dt} = \frac{\partial H}{\partial x}\frac{dx}{dt}+
\frac{\partial H}{\partial y}\frac{dy}{dt}=0.
\end{equation}
It follows that the most general decomposition of this equation is
\begin{equation}\label{eq1.12}
\frac{dx}{dt}=\phi (x,y)\frac{\partial H}{\partial y},\quad
\frac{dy}{dt}=-\phi (x,y)\frac{\partial H}{\partial x},
\end{equation}
where $\phi(x,y)$ is an arbitrary function of $x$ and $y$. Thus,
a priori, an equation of motion for a conservative system leads to a 
unique Hamiltonian, up to an additive constant.
But, a Hamiltonian can correspond to an unlimited number of 
equations of motion.

These results can be generalized to $2N-(x,y)$ variables.

\centerline{\hbox to 3in{\hrulefill}}

Let us now return to \eqn\eqref{eq1.2} and rewrite it to the form
\begin{equation}\label{eq1.13}
\left(\frac23\right)y^2+x^{4/3}=1.
\end{equation}
If we treat this as a functional equation \cite{1}, then the 
following question can be asked: What are possible solutions to
\eqn\eqref{eq1.13}?  Clearly, these possibilities include the following
four piece-wise continuous functions:
\begin{equation}\label{eq1.14}
x(t)=1,\quad y(t)=0,\quad\mbox{for } -\infty < t < +\infty;
\end{equation}
\begin{equation}\label{eq1.15}
x(t)=0,\quad y(t)=\sqrt{\frac32},\quad\mbox{for } -\infty < t < +\infty;
\end{equation}
\begin{equation}\label{eq1.16}
x(t)=\begin{cases}1, & \text{$t>0$},\\
0, & \text{$t<0$},\end{cases},\quad
y(t)=\begin{cases}0, & \text{$t>0$},\\
-\sqrt{\frac32}, & \text{$t<0$};\end{cases}
\end{equation}
and
\begin{subequations}\label{eq1.17}
\begin{equation}
x(t)= \begin{cases}1, & \text{$0<t<\frac T2$},\\
0, & \text{$\frac T2<t<T$},\end{cases}\quad
y(t)= \begin{cases}0, & \text{$0<t<\frac T2$},\\
\sqrt{\frac32}, & \text{$\frac T2<t<T$},\end{cases}\label{eq1.17a}
\end{equation}
{where for $T$ positive}
\begin{equation}
x(t+T)=x(t),\quad y(t+T)=y(t).\label{eq1.17b}
\end{equation}
\end{subequations}

Our purpose, in this paper, is not to investigate the broad variety of 
solutions to \eqn\eqref{eq1.13}, which will be called the 
{\it Leah-functional equation}.  The goal will be to investigate a 
particular class of periodic solutions to this equation, which
follow from making a certain nonlinear transformation of the dependent
variables $x(t)$ and $y(t)$.  The properties of these new periodic
functions, which we call the {\it Imani-functions}, will mimic in large
part mathematical featurs of both the standard trigonometric cosine
and sine functions \cite{3}, and the Jacobian elliptic functions
\cite{2}.

In section 2, we demonstrate the construction of the Imani functions. 
Section~3 presents a summary of the basic properties of these
periodic functions.  Finally, in section~4, we give a brief summary of
our results and several currently unresolved issues related to the
Imani functions.

\section{Construction of the Imani functions}\label{sec2}
Our task is to determine ``some" of the periodic solutions to the 
Leah-functional equation
\begin{equation}\label{eq2.1}
\left(\frac23\right)y(t)^2 + x(t)^{4/3}=1,
\end{equation}
such that they satisfy the conditions
\begin{align}
\mbox{(a)}\qquad & x(0)=1,\quad y(0)=0; \label{eq2.2}\\
\mbox{(b)}\qquad & x(-t)=x(t),\quad y(-t)=-y(t); \label{eq2.3}\\
\mbox{(c)}\qquad & x(t+T)=x(t),\quad y(t+T)=y(t),\quad
\mbox{for } T>0. \label{eq2.4}
\end{align}

To proceed, make the following nonlinear transformation of the
dependent variables
\begin{equation}\label{eq2.5}
u(t)^2=x(t)^{4/3},\quad v(t)^2=\left(\frac23\right)y(t)^2.
\end{equation}
Note that $u(t)$ and $v(t)$ are not uniquely characterized by their 
definitions given in \eqn\eqref{eq2.5}; there are $(\pm)$ signs which can 
appear.  However, this ambiquity will be used to enforce the condition
specified in \eqn\eqref{eq2.3}.

From \eqns\eqref{eq2.1} and \eqref{eq2.5}, it follows that
\begin{equation}\label{eq2.6}
u(t)^2+v(t)^2=1.
\end{equation}
Now make the following identification
\begin{equation}\label{eq2.7}
u(t)^2=[\cos \psi (t)]^2,\quad
v(t)^2=[\sin \psi (t)]^2.
\end{equation}
To satisfy the conditions in \eqns\eqref{eq2.3} and \eqref{eq2.4},
$\psi(t)$ is required to have the properties
\begin{subequations}\label{eq2.8}
\begin{align}
\psi (-t)&=-\psi (t)\label{eq2.8a}\\
\psi (t+T)&= \psi (t)+2\pi.\label{eq2.8b}
\end{align}
\end{subequations}
Consequently, $x(t)$ and $y(t)$ may be selected to be 
\begin{subequations}\label{eq2.9}
\begin{align}
x(t)&= \left[\mbox{sgn}(\cos\psi (t))\right]
|\cos\psi (t)|^{3/2}\label{eq2.9a}\\
y(t)&= \sqrt{\frac32}\sin \psi (t).\label{eq2.9b}
\end{align}
\end{subequations}
In the above expressions, ``$|\cdots|$" denotes the absolute value
and ``sgn" is the sign-function \cite{6}, i.e.,
\begin{equation}\label{eq2.10}
\mbox{sgn} (t)\equiv \begin{cases}
+1, & \text{$t>0$};\\
0, & \text{$t=0$};\\
-1, & \text{$t<0$}.\end{cases}
\end{equation}

Observe that $\psi (t)$ satisfies the functional equation given in
\eqn\eqref{eq2.8} and its solutions is \cite{6,8}
\begin{subequations}\label{eq2.11}
\begin{equation}
\psi (t)=A(t)+\left(\frac{2\pi}T\right)t,\label{eq2.11a}
\end{equation}
where
\begin{equation}
A(-t)=-A(t),\quad A(t+T)=A(t).\label{eq2.11b}
\end{equation}
\end{subequations}
Assuming $A(t)$ has a Fourier series, it follows that \cite{8}
\begin{equation}\label{eq2.12}
A(t)=\sum^\infty_{k=1}a_k\sin\left(\frac{2\pi k}T\right)t.
\end{equation}

We name the functions given in \eqns\eqref{eq2.9},
respectively, the {\it Imani cosine and sine functions}, and
denote them by the symbols ``Ics" and ``Isn".

\section{Properties of the Imani periodic functions}\label{sec3}

The Imani functions are a class of continuous, periodic solutions to the 
Leah functional equation.  By construction, they satisfy the conditions
specified in \eqns\eqref{eq2.2}, \eqref{eq2.3}, and \eqref{eq2.4}. 
Note that these functions may have a period $T$, of any magnitude,
for $T>0$.  Central to these representations, for Ics$(t)$ and
Isn$(t)$, is the fact that $\psi (t)$ has the mathematical form given in
\eqns\eqref{eq2.11a} and \eqref{eq2.12}.

Overall, except for the complexity of the final representation
for Ics$(t)$, the general properties of the Imani functions are
similar to the standard trigonometric cosine and sine functions.  
It should also be clear that the construction of this paper could
have been done using the Jacobi cosine and sine functions \cite{2}.

\begin{figure}[H]
\centerline{\includegraphics[width=4.50in,height=2.8in]{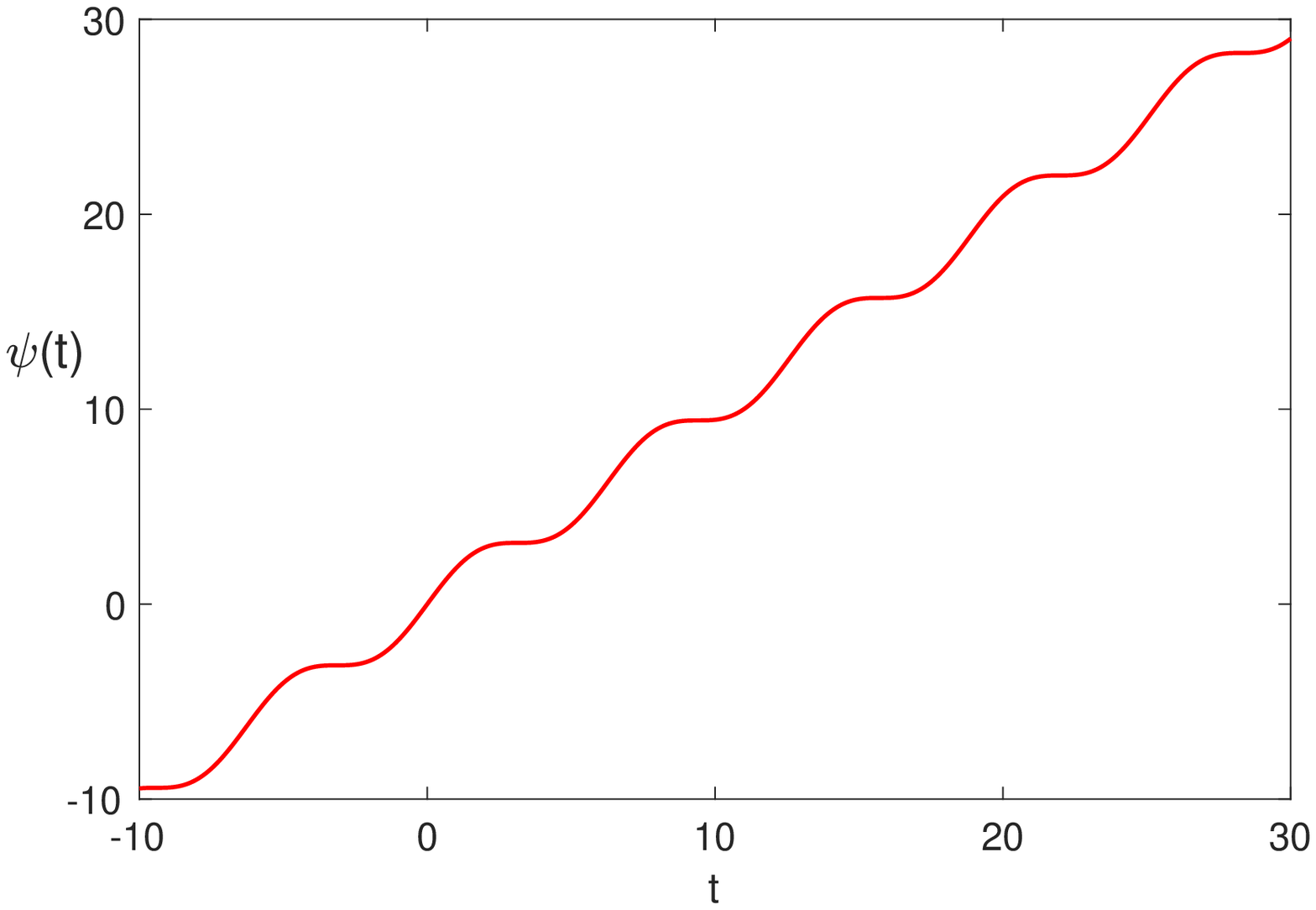}}
\caption{}
\label{fig1}
\end{figure}

\begin{figure}[H]
\centerline{\includegraphics[width=4.50in,height=3in]{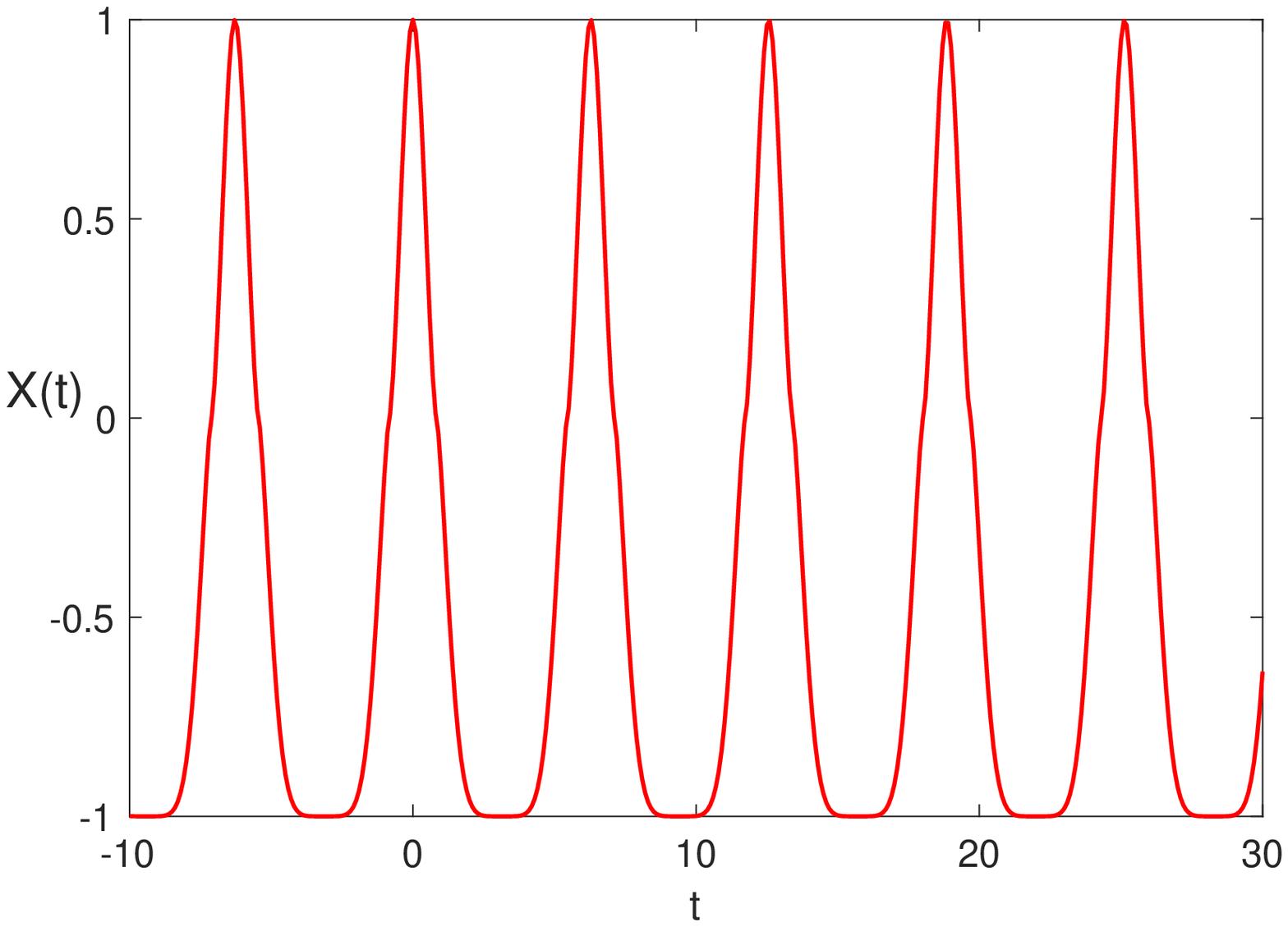}}
\caption{}
\label{fig2}
\end{figure}

Figure~\ref{fig1} is a plot
of $\psi (t)$ vs $t$ , while Figure~\ref{fig2} 
is the corresponding plot of $x(t)$ vs $t$.
Both of these results are for a single term in the Fourier expansion of 
$A(t)$; see Equation~\eqref{eq2.12}. 

\section{Future work}\label{sec4}

In summary, we have shown how to construct a particular set of periodic
solutions to the Leah functional equation.  The following questions
can be asked:

1) What other classes of continuous, periodic functions exist?

2) What is the full Fourier sereies representation for the 
Imani cosine function?

3) Do derivatives exist for Ics$(t)$ and Isn$(t)$? 
If so, what are they?

4) What are the integrals
\begin{equation}\label{eq4.1}
\int\mbox{Ics}(t) dt\quad\mbox{and}\quad \int\mbox{Isn}(t)dt
\end{equation}
of these functions?

5) Can products of Imani functions
\begin{equation}\label{eq4.2}
\mbox{Ics}(t_1)\mbox{Ics}(t_2),\quad
\mbox{Isn}(t_1)\mbox{Isn}(t_2),\quad
\mbox{Isn}(t_1)\mbox{Ics}(t_2),
\end{equation}
be calculated?

6) What are Ics$(t_1+t_2)$ and Isn$(t_1+t_2)$?

7) Are the Imani functions solutions of differential equations?
If so, what are they? Is there a single differential equation for 
which both are special solutions?

8) Finally, can the Imani periodic functions be analyzed using
geometrical methods such as those used for the elliptic functions
\cite{11} and the square periodic functions \cite{7}?

\section*{Acknowledgments}

We thank Professors Torini Lewis and Sandra A. Rucker, both at
Clark Atlanta University, Department of Mathematical Sciences, and
Kale Oyedeji, Department of Physics, Morehouse College, for several
insightful and productive discussions on periodic functions.

\end{document}